\documentclass[11pt, twoside, leqno]{article}

\usepackage{amsmath,amsthm}
\usepackage{amssymb}

\usepackage{graphicx}

\usepackage[T1]{fontenc}

\usepackage[margin=1in]{geometry}
\usepackage{latexsym, amsfonts, amsmath,tikz}
\usepgflibrary{arrows}
\newcommand{\be}{\begin{equation}}
\newcommand{\ee}{\end{equation}}
\newcommand{\bea}{\begin{eqnarray}}
\newcommand{\eea}{\end{eqnarray}}
\newcommand{\bean}{\begin{eqnarray*}} 
\newcommand{\eean}{\end{eqnarray*}}

\newcommand{\brray}{\begin{array}}
\newcommand{\erray}{\end{array}}
\newcommand{\ben}{\begin{equation}{nonumber}}
\newcommand{\een}{\end{equation}{nonumber}}

\newtheorem{dfn}{Definition}[section]
\newtheorem{thm}[dfn]{Theorem}
\newtheorem{lmma}[dfn]{Lemma}
\newtheorem{ppsn}[dfn]{Proposition}
\newtheorem{crlre}[dfn]{Corollary}
\newtheorem{xmpl}[dfn]{Example}
\newtheorem{rmrk}[dfn]{Remark}
\newcommand{\bdfn}{\begin{dfn}}
\newcommand{\bthm}{\begin{thm}}
\newcommand{\blmma}{\begin{lmma}}
\newcommand{\bppsn}{\begin{ppsn}}
\newcommand{\bcrlre}{\begin{crlre}}
\newcommand{\bxmpl}{\begin{xmpl}}
\newcommand{\brmrk}{\begin{rmrk}}
\newcommand{\edfn}{\end{dfn}}
\newcommand{\ethm}{\end{thm}}
\newcommand{\elmma}{\end{lmma}}
\newcommand{\eppsn}{\end{ppsn}}
\newcommand{\ecrlre}{\end{crlre}}
\newcommand{\exmpl}{\end{xmpl}}
\newcommand{\ermrk}{\end{rmrk}}

\newcommand{\cla}{{\cal A}}
\newcommand{\clb}{{\cal B}}
\newcommand{\clc}{{\cal C}}
\newcommand{\cld}{{\cal D}}
\newcommand{\cle}{{\cal E}}
\newcommand{\clf}{{\cal F}}

\newcommand{\clh}{{\cal H}}
\newcommand{\cli}{{\cal I}}
\newcommand{\clj}{{\cal J}}

\newcommand{\cll}{{\cal L}}

\newcommand{\cln}{{\cal N}}
\newcommand{\clo}{{\cal O}}
\newcommand{\clp}{{\cal P}}
\newcommand{\clq}{{\cal Q}}

\newcommand{\cls}{{\cal S}}
\newcommand{\clt}{{\cal T}}

\newcommand{\clv}{{\cal V}}
\newcommand{\clw}{{\cal W}}

\newcommand{\clz}{{\cal Z}}

\def\a*{{\cal A}_{h,*}}
\def\B{{\cal B}(h)}
\def\B1{{\cal B}_1(h)}
\def\b{{\cal B}^{\rm s.a.}(h)}
\def\b1{{\cal B}^{\rm s.a.}_1(h)}

\newcommand{\ot}{\otimes}

\newcommand{\raro}{\rightarrow}

\def \qed {$\Box$}
\mathchardef\mhyphen="2D 

\newcommand{\midarrow}{\tikz \draw[-triangle 90] (0,0) -- +(.1,0);}

\def\a*{{\cal A}_{h,*}}
\def\B{{\cal B}(h)}
\def\B1{{\cal B}_1(h)}
\def\b{{\cal B}^{\rm s.a.}(h)}
\def\b1{{\cal B}^{\rm s.a.}_1(h)}

\mathchardef\mhyphen="2D 
\begin{document}

	\baselineskip=17pt
	

	\title{Improved bound on the non zero eigenvalues of the graph Laplacian coming from the quantum symmetry of vertex transitive graphs}
	
	\author{Soumalya Joardar\\
	IISER Kokata, Mohanpur\\
West Bengal-741246, India\\
email: soumalya.j@gmail.com}
	
	\date{}
	
	\maketitle
	\renewcommand{\thefootnote}{}
	\footnote{2010 \emph{Mathematics Subject Classification}: 46L30; 05C50.}
	\footnote{\emph{Key words and phrases}: Compact quantum group, graph Laplacian, quantum symmetry}
	\footnote{The author acknowledges the support from the Department of Science and Technology, India (DST/INSPIRE/04/2016/002469)}
	\renewcommand{\thefootnote}{arabic{footnote}}
	\setcounter{footnote}{0}
	\begin{abstract}
	A chain of quantum subgroups of the quantum automorphism group of finite graphs has been introduced. It generalizes the construction of J. Bichon (see \cite{bichon}) in a sense. A better bound of the non zero eigenvalues of the graph Laplacian has been obtained using the chain of quantum subgroups. 
	\end{abstract}
\section{Introduction}
The role of compact quantum groups (CQG in short) as the symmetry object in the realm of noncommutative geometry is now well established. Study of quantum symmetry of various mathematical structures potentially has two applications. In one hand it can enrich the theory of CQG's. On the other hand, it can supply more information of the mathematical structure in question than what its classical symmetry does. Some of the interesting mathematical structures come from finite graphs. The study of quantum symmetry of finite structures including graphs was initiated by S. Wang, J. Bichon, T. Banica (see \cite{wang}, \cite{bichon}, \cite{banica}) among others. They were mainly interested in studying the quantum group theoretic aspects of such symmetries. In fact it is perhaps fair to say that until now almost all the work in quantum symmetry of graphs has found its use in understanding the theory of compact quantum groups and their representation (see also \cite{bancolo}, \cite{Maccinska}). Given this context, this article tries to use quantum symmetry to extract spectral information of finite, connected graphs. \\
\indent It is well known that the graph Laplacian (normalized in the sense of \cite{Chung}) has non-negative eigenvalues. The lower bound of the positive eigenvalues provides a great deal of information about the connectivity of the graph. For a symmetric graph (classical vertex transitive), using the symmetry of the graph one can obtain lower bound of the positive eigenvalues. In this article  using the quantum symmetry of a vertex transitive, finite, connected graph we improve the bound of the positive eigenvalues of the graph Laplacian. We adapt the techniques from \cite{Chung} to the quantum set up to obtain our main result (Theorem \ref{main}).\\
\indent On the quantum group theoretic front, it is known that the action of the quantum automorphism group of finite graphs in the sense of Banica is $2$-transitive on the set of vertices, but it fails to be higher order transitive in general (see \cite{Maccinska}, \cite{bancolo}). This can be attributed to the fact that unlike the classical automorphism group, the quantum automorphism group (in the sense of Banica) does not act on the path spaces of any length. We introduce a chain of quantum subgroups of the quantum automorphism group in the sense of Banica containing the classical automorphism group. It generalizes the construction of Bichon (in \cite{bichon}) in a sense. The members of the chain act transitively on the path spaces of certain length.\\    
\indent Now let us discuss the organization of the paper briefly. The paper is divided into two sections viz. the preliminary section and the main section. In the preliminary section we mainly collect some well known facts about the quantum symmetry of finite graphs as well as the eigenvalues of the Laplacian of finite graphs. In particular the classical result on the lower bound of the eigenvalues of the Laplacian of vertex transitive graph using the classical automorphism group is recalled. The main section is divided into two subsections. In one subsection, a chain of quantum subgroups containing the classical automorphism group is introduced and the question of higher order path transitivity has been considered. The last subsection is devoted to obtaining a better lower bound on the non zero eigenvalues of the graph Laplacian of vertex transitive graphs adapting the techniques of \cite{Chung} to the quantum set up.  
\section{Preliminaries}
\subsection{Eigenvalues of the graph Laplacian}
\label{prel_graph}
A finite graph $\Gamma$ is a collection of finitely many vertices and edges. We shall denote the set of vertices by $V$ and the set of edges by $E$. The vertices will be denoted by the natural numbers. Given a finite graph $\Gamma$ with $n$-vertices and without multiple edges, the vertex matrix is given by a $0-1$-valued $n\times n$ matrix $\epsilon=((\epsilon_{ij}))$ such that the $ij$-th entry is $1$ if there is an edge between the vertices $i$ and $j$, zero otherwise. In this paper, all the graphs are without loops. For a vertex $i\in V$, degree of the vertex is given by the number of edges emitting from the vertex $i$ and is denoted by $d_{i}$. Volume of the graph is defined to be
\begin{displaymath}
\clv=\sum_{i\in V}d_{i}.
\end{displaymath}
Given a graph $\Gamma$, consider all the shortest paths between any two vertices. The diameter is defined to be the maximum of the number of edges in the shortest paths. It is denoted by $D$. For any two vertices $i,j$, we write $i\sim j$ if $\epsilon_{ij}=1$. If the symmetry group of a graph $\Gamma$ is denoted by ${\rm Aut}(\Gamma)$, then a graph is said to be vertex transitive if for any two vertices $i,j$ there is a permutation in ${\rm Aut}(\Gamma)$ which takes $i$ to $j$. The automorphism group of the graph defines an equivalence relation on the set of edges. Two edges $e_{i}, e_{j}$ are said to be equivalent if there is an element of the automorphism group which takes $e_{i}$ to $e_{j}$. Then one can consider equivalence classes of edges denoted by $E_{1},...,E_{m}$ (say). The index of a graph is defined to be
\begin{displaymath}
{\rm ind}(\Gamma)=\frac{\clv}{2{\rm min}_{i}|E_{i}|}
\end{displaymath} Now we recall the Laplacian (normalized in the sense of \cite{Chung}) on a finite, connected graph without loops and multiple edges. The Laplacian is an $n\times n$ martrix given by
\begin{eqnarray*}
	\cll(i,j)=\begin{cases}
		1, \ if \ i=j \ and \ d_{i}\neq 0\\
		-\frac{1}{\sqrt{d_{i}d_{j}}}, \ i\sim j\\
		0, \ else.
	\end{cases} 
\end{eqnarray*}
It can be shown that the eigenvalues are non negative. The eigen values are written as $0=\lambda_{0}\leq\lambda_{1}\leq...\leq\lambda_{n-1}$. For a vertex transitive graph $\Gamma$ with diameter $D$, we have the following lower bound of the first non-zero eigenvalue (see Theorem 7.7 of \cite{Chung}):
\begin{eqnarray}
\label{bound}
\lambda_{1}\geq\frac{1}{D^{2}{\rm ind}(\Gamma)}. 
\end{eqnarray}
\subsection{Quantum automorphism group of finite graphs}
\bdfn (see \cite{Woro}, \cite{Van})
A compact quantum group (CQG for short) $\mathbb{G}$ is a pair $(C(\mathbb{G}),\Delta)$ such that $C(\mathbb{G})$ is a unital $C^{\ast}$-algebra and $\Delta:C(\mathbb{G})\raro C(\mathbb{G})\ot_{\rm min}C(\mathbb{G})$ is a unital $C^{\ast}$-homomorphism such that\\
(i) $({\rm id}\ot\Delta)\circ\Delta=(\Delta\ot{\rm id})\circ\Delta.$\\
(ii) Linear span of the sets $\Delta(C(\mathbb{G}))(1\ot C(\mathbb{G}))$ and $\Delta(C(\mathbb{G}))(C(\mathbb{G})\ot 1)$ are norm dense in $C(\mathbb{G})\ot_{\rm min}C(\mathbb{G})$.
\edfn
Let $W$ be a finite dimensional vector space. Given a $\mathbb{C}$-linear map $u:W\raro W\ot_{\mathbb{C}}C(\mathbb{G})$ recall the leg numbering notations $u_{12}$ and $u_{13}$.
\bdfn
A corepresentation of a CQG $\mathbb{G}$ on a finite dimensional vector space $W$ is a $\mathbb{C}$-linear map $u:W\raro W\ot_{\mathbb{C}}C(\mathbb{G})$ such that $({\rm id}\ot \Delta)\circ u=u_{12}u_{13}$. 
\edfn
Let $\{e_{i}\}_{i=1,...,n}$ be  a basis of $W$. Let us write $u(e_{i})=\sum_{j=1}^{n}e_{j}\ot q_{ji}$ for $q_{ji}\in C(\mathbb{G})$. Then $u$ is said to be non degenerate if the matrix $((q_{ij}))\in M_{n}(C(\mathbb{G}))$ is invertible. Given a finite dimensional non degenerate representation $u$ of $\mathbb{G}$ on a finite dimensional vector space $W$, we have the corresponding non degenerate $\gamma$-th tensor product representation to be denoted by $u^{\ot\gamma}$ on the finite dimensional vector space $W\underbrace{\ot...\ot}_{\gamma \ {\rm times}}W$. The $\mathbb{C}$-linear map $u^{\ot\gamma}$ is given on the basis elements by\begin{displaymath}
u^{\ot\gamma}(e_{i_{1}}\ot...\ot e_{i_{\gamma}})=\sum_{j_{1},...,j_{\gamma}}e_{j_{1}}\ot...\ot e_{j_{\gamma}}\ot q_{j_{1}i_{1}}...q_{j_{\gamma}i_{\gamma}}.
\end{displaymath} 
\bdfn
A CQG $\mathbb{G}$ is said to coact on a unital $C^{\ast}$-algebra $\mathbb{A}$ if there is a unital $C^{\ast}$-homomorphism $u:\mathbb{A}\raro \mathbb{A}\ot_{\rm min}C(\mathbb{G})$ such that\\
(i) $({\rm id}\ot\Delta)\circ u=(u\ot {\rm id})\circ u$.\\
(ii) Linear span of the set $u(\mathbb{A})(1\ot C(\mathbb{G}))$ is norm dense in $\mathbb{A}\ot_{\rm min}C(\mathbb{G})$.
\edfn 
\bxmpl
The quantum permutation group is the `free' version of the classical permutation group. It is denoted by $S_{n}^{+}$. The underlying $C^{\ast}$-algebra $C(S_{n}^{+})$ which is infinite dimensional for $n\geq 4$, is the universal $C^{\ast}$-algebra generated by selfadjoint elements $q_{ij}$'s such that
\begin{eqnarray}
\sum_{j}q_{ij}=1= \sum_{j}q_{ji} \ \forall i, q_{ij}q_{ik}=\delta_{j,k}q_{ij} \ \forall i, q_{ji}q_{ki}=\delta_{j,k}q_{ji} \ \forall i.
\end{eqnarray}
The coproduct $\Delta$ is given on the generators by $\Delta(q_{ij})=\sum_{k=1}^{n}q_{ik}\ot q_{kj}$. The CQG coacts on the finite dimensional $C^{\ast}$-algebra $C(X)$ where $X$ is a set with $n$-points. The corresponding coaction is given by 
\begin{eqnarray}
u(\chi_{i})=\sum_{j=1}^{n}\chi_{j}\ot q_{ji}.
\end{eqnarray} In fact $S_{n}^{+}$ is the universal object in the category of CQG's coacting on the $C^{\ast}$-algebra $C(X)$ (see \cite{wang}). The $\mathbb{C}$-linear map $u$ can also be viewed as a finite dimensional corepresentation of $S_{n}^{+}$ on the $n$-dimensional vector space $C(X)$. Note that $\chi_{i}$ denotes the function which takes the value $1$ on the point $i$ and zero elsewhere. It is easy to see that $\{\chi_{i}\}_{i=1,...,n}$ is a linear basis of $C(X)$.
\exmpl
Recall the definition of a finite connected graph $\Gamma$ from the Subsection \ref{prel_graph}. 
\bdfn (See \cite{banica})
For a finite graph $\Gamma$ without multiple edges, the quantum automorphism group of $\Gamma$ in the sense of Banica, to be denoted by ${\rm Aut}^{+}(\Gamma)$, is the CQG whose underlying $C^{\ast}$-algebra $C({\rm Aut}^{+}(\Gamma))$ is given by the $C^{\ast}$-algebra $\frac{C(S_{n}^{+})}{<q\epsilon-\epsilon q>}$, where $q\in M_{n}(C(S^{+}_{n}))$ is the matrix $((q_{ij}))$. The coaction of ${\rm Aut}^{+}(\Gamma)$ on $C(V)$ is again given by
\begin{displaymath}
u(\chi_{i})=\sum_{j=1}^{n}\chi_{j}\ot q_{ji},
\end{displaymath}
where $\chi_{i}$ is the function which takes the value $1$ on the vertex $i$ and zero on the other vertices.
\edfn
For the description of the $C^{\ast}$-algebra $C({\rm Aut}^{+}(\Gamma))$ in terms of generators and relations, see \cite{Schmidt}. We shall use the following relation among the generators in particular in this paper:\\
For $(i_{1},i_{2})\in E$ and $(j_{1},j_{2})\notin E$, we have \begin{eqnarray}
\label{Banica}q_{i_{1}j_{1}}q_{i_{2}j_{2}}=q_{i_{2}j_{2}}q_{i_{1}j_{1}}=q_{j_{1}i_{1}}q_{j_{2}i_{2}}=q_{j_{2}i_{2}}q_{j_{1}i_{1}}=0.\end{eqnarray}
 Let us denote the set of paths of length $\alpha$ by $E^{\alpha}$. Then the space $C(E^{\alpha})$ is a linear subspace of the space $C(V)\underbrace{\ot.........\ot}_{(\alpha+1) \ {\rm times}}C(V)\cong C(V\underbrace{\times.........\times}_{(\alpha+1) \ {\rm times}}V)$. More precisely $C(E^{\alpha})={\rm Sp}\{\chi_{i_{1}}\ot...\ot\chi_{i_{\alpha+1}}:(i_{1},...,i_{\alpha+1}) \ {\rm is \ a \ path}\}$. We denote the basis element $\chi_{i_{1}}\ot...\ot\chi_{i_{\alpha+1}}$ by $\chi_{i_{1}...i_{\alpha+1}}$. In fact we have the following direct sum decomposition of the vector space $W=C(V)\ot...\ot C(V)$:
 \begin{displaymath}
 C(V)\ot...\ot C(V)=C(E^{\alpha})\oplus W_{0}, 
 \end{displaymath}
 where $W_{0}={\rm Sp}\{\chi_{i_{1}}\ot...\ot\chi_{i_{\alpha+1}}:(i_{1},...,i_{\alpha+1}) {\rm is \ not \ a \ path}\}$. Now as before considering $u$ as a nondegenerate corepresentation of ${\rm Aut}^{+}(\Gamma)$ on $C(V)$, we have the corresponding corepresentation $u^{\ot(\alpha+1)}$.
 \blmma
 \label{reponpath}
 The corepresentation $u^{\ot(\alpha+1)}$ restricts to a nondegenerate corepresentation on $C(E^{\alpha})$ of the CQG ${\rm Aut}^{+}(\Gamma)$.
 \elmma 
 \begin{proof}
 	Let $(i_{1},...,i_{\alpha+1})$ be a path. Then by definition of the corepresentation $u^{\ot(\alpha+1)}$, we have
 	\begin{displaymath}
 	u^{\ot(\alpha+1)}(\chi_{i_{1}}\ot...\ot\chi_{i_{\alpha+1}})=\sum_{j_{1},...,j_{\alpha+1}}\chi_{j_{1}}\ot...\ot\chi_{j_{\alpha+1}}\ot q_{j_{1}i_{1}}...q_{j_{\alpha+1}i_{\alpha+1}}.
 	\end{displaymath}
 	But by equation (\ref{Banica}), for the indices $(j_{1},...,j_{\alpha+1})$ such that $(j_{1},...,j_{\alpha+1})$ are not paths, \begin{displaymath}q_{j_{1}i_{1}}...q_{j_{\alpha+1}i_{\alpha+1}}=0.\end{displaymath} This implies that $u^{\ot(\alpha+1)}(C(E^{\alpha}))\subset C(E^{\alpha})\ot C({\rm Aut}^{+}(\Gamma))$. By similar reasoning we can show that $u^{\ot(\alpha+1)}(W_{0})\subset W_{0}\ot C({\rm Aut}^{+}(\Gamma))$ proving that in deed $C(E^{\alpha})$ is a reducing subspace.
 	\end{proof}

\section{Main section}   

\subsection{A chain of quantum subgroups}
Given a finite graph without any multiple edge, the classical automorphism group $\mathrm{Aut}(\Gamma)$ is a quantum subgroup of ${\rm Aut}^{+}(\Gamma)$. Unlike the classical automorphism group, ${\rm Aut}^{+}(\Gamma)$ does not act on the path spaces (although it has a corepresentation on the path space as observed) of the graph. We are going to introduce a CQG for a natural number $k$ such that it coacts on the space of paths of length less than or equal to $k$. For $k\in\mathbb{N}$, let $C(\mathbb{G}_{k}^{+}(\Gamma))$ be the $C^{\ast}$-algebra generated by $\{q_{ij}:i,j=1,...,n\}$ such that $q_{ij}$'s satisfies the relations of $C({\rm Aut}^{+}(\Gamma))$ along with the following relations:
\begin{eqnarray}
\label{1}
q_{i_{1}j_{1}}...q_{i_{\beta}j_{\beta}}=q_{i_{\beta}j_{\beta}}...q_{i_{1}j_{1}}, \beta\leq (k+1)
\end{eqnarray}
\blmma
For a finite graph $\Gamma$ without any multiple edge, $(C(\mathbb{G}_{k}^{+}(\Gamma)),\Delta)$ is a quantum subgroup of ${\rm Aut}^{+}(\Gamma)$ for all $k\in\mathbb{N}$, where the coproduct $\Delta$ is again given on the generators by $\Delta(q_{ij})=\sum_{k=1}^{n}q_{ik}\ot q_{kj}$.
\elmma
\begin{proof}
	Clearly it suffices to show that $\Delta(q_{i_{1}j_{1}}...q_{i_{\beta}j_{\beta}}-q_{i_{\beta}j_{\beta}}...q_{i_{1}j_{1}})=0$ for all $\beta\leq (k+1)$. So fix a $\beta\leq(k+1)$.
	\begin{eqnarray*}
	&&\Delta(q_{i_{1}j_{1}}...q_{i_{\beta}j_{\beta}}-q_{i_{\beta}j_{\beta}}...q_{i_{1}j_{1}})\\
		&=& \sum_{(k_{1},...,k_{\beta}) \ path}q_{i_{1}k_{1}}...q_{i_{\beta}k_{\beta}}\ot q_{k_{1}j_{1}}...q_{k_{\beta}j_{\beta}}-\sum_{(l_{1},...,l_{\beta})}q_{i_{\beta}l_{1}}...q_{i_{1}l_{\beta}}\ot q_{l_{1}j_{\beta}}...q_{l_{\beta}j_{1}}.
\end{eqnarray*}
Observe that $q_{i_{\beta}l_{1}}...q_{i_{1}l_{\beta}}$ is zero if $(l_{\beta},...,l_{1})$ is not a path. Hence we get \begin{displaymath}\Delta(q_{i_{1}j_{1}}...q_{i_{\beta}j_{\beta}}-q_{i_{\beta}j_{\beta}}...q_{i_{1}j_{1}})=0.\end{displaymath}
		\end{proof}
 It is straightforward to observe the following inclusions:
   \begin{displaymath}
   {\rm Aut}(\Gamma)\leq...\leq\mathbb{G}_{k+1}^{+}(\Gamma)\leq\mathbb{G}_{k}^{+}(\Gamma)\leq...\leq\mathbb{G}_{1}^{+}(\Gamma)\leq{\rm Aut}^{+}(\Gamma).
   \end{displaymath}
   Note that $\mathbb{G}^{+}_{1}(\Gamma)$ is the quantum automorphism group of a finite graph $\Gamma$ without any multiple edge in the sense of Bichon (see \cite{bichon}). We shall denote the image of $q_{ij}\in{\rm Aut}^{+}(\Gamma)$ in any of the $\mathbb{G}_{k}^{+}(\Gamma)$'s again by $q_{ij}$.
  \blmma
\label{key}
Fix a $k\in\mathbb{N}$ and $\beta\leq (k+1)$. Then for $q_{ij}$'s in $\mathbb{G}_{k}^{+}(\Gamma)$, we have the following:
\begin{displaymath}
q_{i_{1}j_{1}}...q_{i_{\beta}j_{\beta}}q_{i_{1}l_{1}}...q_{i_{\beta}l_{\beta}}=\delta_{j_{1}l_{1}}...\delta_{j_{\beta}l_{\beta}}q_{i_{1}j_{1}}...q_{i_{\beta}j_{\beta}},
\end{displaymath}
where $(i_{1},...,i_{\beta}), (j_{1},...,j_{\beta}), (l_{1},...,l_{\beta})$ are paths.
\elmma
\begin{proof}
	Since all the $q_{ij}$'s are in $C(\mathbb{G}^{+}_{k}(\Gamma))$, we get
	\begin{eqnarray*}
	q_{i_{1}j_{1}}...q_{i_{\beta}j_{\beta}}q_{i_{1}l_{1}}...q_{i_{\beta}l_{\beta}}&=&q_{i_{1}j_{1}}...q_{i_{\beta}j_{\beta}}q_{i_{\beta}l_{\beta}}...q_{i_{1}l_{1}}\\
	&=&\delta_{j_{\beta}l_{\beta}}q_{i_{1}j_{1}}...q_{i_{\beta-1}j_{\beta-1}}q_{i_{\beta}j_{\beta}}q_{i_{\beta-1}l_{\beta-1}}...q_{i_{1}l_{1}}\\
	&=& \delta_{j_{\beta}l_{\beta}}q_{i_{1}j_{1}}...q_{i_{\beta}j_{\beta}}q_{i_{\beta-1}j_{\beta-1}}q_{i_{\beta-1}l_{\beta-1}}...q_{i_{1}l_{1}}\\
	&=& \delta_{j_{\beta}l_{\beta}}\delta_{j_{\beta-1}l_{\beta-1}}q_{i_{1}j_{1}}...q_{i_{\beta}j_{\beta}}q_{i_{\beta-1}j_{\beta-1}}q_{i_{\beta-2}l_{\beta-2}}...q_{i_{1}l_{1}}\\
	&=& \delta_{j_{\beta}l_{\beta}}\delta_{j_{\beta-1}l_{\beta-1}}q_{i_{1}j_{1}}...q_{i_{\beta-2}j_{\beta-2}}q_{i_{\beta-1}j_{\beta-1}}q_{i_{\beta}j_{\beta}}q_{i_{\beta-2}l_{\beta-2}}...q_{i_{1}l_{1}}
	\end{eqnarray*}
Continuing thus and repetitively using (\ref{1}), we finish the proof of the lemma.
\end{proof}
\blmma
\label{action}
$\mathbb{G}^{+}_{k}(\Gamma)$ acts on $C(E^{\gamma})$ for $\gamma\leq k$.
\elmma 
\begin{proof}
	$C(E^{\gamma})$ is generated by self adjoint mutually orthogonal projections $\chi_{i_{1}...i_{\gamma+1}}$ such that $(i_{1},...,i_{\gamma+1})$ is a path and $\sum_{(i_{1},...,i_{\gamma+1})}\chi_{i_{1}...i_{\gamma+1}}=1$. Recall the corepresentation $u^{\ot(\gamma+1)}$ of $\mathbb{G}_{k}^{+}$ (being a quantum subgroup of ${\rm Aut}^{+}(\Gamma)$ for all $k$) on $C(E^{\gamma})$. We shall prove that $u^{\ot(\gamma+1)}$ is actually a $C^{\ast}$-coaction. To that end,
	\begin{displaymath}\sum_{(i_{1},...,i_{\gamma+1}) \ path}u^{\ot(\gamma+1)}(\chi_{i_{1}...i_{\gamma+1}})=\sum_{(j_{1},...,j_{\gamma+1}) \ path}\chi_{j_{1}...j_{\gamma+1}}\ot \sum_{(i_{1},...,i_{\gamma+1}) \ path}q_{j_{1}i_{i}}...q_{j_{\gamma+1}i_{\gamma+1}}.\end{displaymath} We can replace the above summation over $(i_{1},...,i_{\gamma+1})$ such that $(i_{1},...,i_{\gamma+1})$ path by all tuples $(i_{1},...,i_{\gamma+1})$ since for $(i_{1},...,i_{\gamma+1})$ not a path $q_{j_{1}i_{1}}...q_{j_{\gamma+1}i_{\gamma+1}}=0$. Using this and the magic unitarity of $q_{ij}$, we get $\sum_{(i_{1},...,i_{\gamma+1}) \ path}u^{\ot(\gamma+1)}(\chi_{i_{1}...i_{\gamma+1}})=1\ot 1$. Using (\ref{1}), we get $u^{\ot(\gamma+1)}(\chi_{i_{1}...i_{\gamma+1}})^{\ast}=u^{\ot(\gamma+1)}(\chi_{i_{1}...i_{\gamma+1}})$. Using the mutual orthogonality of $\chi_{i_{1}....i_{\gamma+1}}$, we get
	\begin{eqnarray*}
		&& u^{\ot(\gamma+1)}(\chi_{i_{1}...i_{\gamma+1}})u^{\ot(\gamma+1)}(\chi_{j_{1}...j_{\gamma+1}})\\
		&=& \sum_{(l_{1},...,l_{\gamma+1}) \ path} \chi_{l_{1}...l_{\gamma+1}}\ot q_{l_{1}i_{1}}...q_{l_{\gamma+1}i_{\gamma+1}}q_{l_{1}j_{1}}...q_{l_{\gamma+1}j_{\gamma+1}}\\
		&=& \delta_{i_{1}j_{1}}...\delta_{i_{\gamma+1}j_{\gamma+1}}\sum_{(l_{1},...,l_{\gamma+1}) \ path} \chi_{l_{1}...l_{\gamma+1}}\ot q_{l_{1}i_{1}}...q_{l_{\gamma+1}i_{\gamma+1}} (by \ Lemma \ \ref{key})\\
		&=& \delta_{i_{1}j_{1}}...\delta_{i_{\gamma+1}j_{\gamma+1}}u^{\ot(\gamma+1)}(\chi_{i_{1}...i_{\gamma+1}}).
		\end{eqnarray*} 
This proves in deed $u^{\ot(\gamma+1)}$ is a well defined $C^{\ast}$-homomorphism. Coassociativity and density condition follow from the same of the corresponding property of representation.
\end{proof}
\indent Now we define a relation on the path space by the following:\\
Fix a $k\in\mathbb{N}$. Then two paths $(i_{1},...,i_{\gamma})$ and $(j_{1},...,j_{\gamma})$ are said to be $k$-equivalent and denoted by $(i_{1},...,i_{\gamma})\sim^{k}(j_{1},...,j_{\gamma})$ if
$q_{i_{1}j_{1}}...q_{i_{\gamma}j_{\gamma}}\neq 0$ for $q_{ij}\in C(\mathbb{G}^{+}_{k}(\Gamma))$.
\blmma
$\sim^{k}$ is an equivalence relation on the space of paths of length less than or equal to $(2k+1)$.
\elmma 
\begin{proof}
	Symmetry and reflexivity are clear (see \cite{Maccinska}). The only non trivial part is the transitivity which we show for paths of length $(2\gamma+1)$ for $\gamma\leq k$. To that end suppose $(i_{1},...,i_{2\gamma+2})\sim^{k}(j_{1},...,j_{2\gamma+2})\sim^{k}(l_{1},...,l_{2\gamma+2})$.
	\begin{eqnarray*}
		&&(q_{i_{1}j_{1}}...q_{i_{\gamma+1}j_{\gamma+1}}\ot q_{j_{1}l_{1}}...q_{j_{\gamma+1}l_{\gamma+1}})\Delta(q_{i_{1}l_{1}}...q_{i_{\gamma+1}l_{\gamma+1}}...q_{i_{2\gamma+2}l_{2\gamma+2}})\\
		&&(q_{i_{\gamma+2}j_{\gamma+2}}...q_{i_{2\gamma+2}j_{2\gamma+2}}\ot q_{j_{\gamma+2}l_{\gamma+2}}...q_{j_{2\gamma+2}l_{2\gamma+2}})\\
		&=&(q_{i_{1}j_{1}}...q_{i_{\gamma+1}j_{\gamma+1}}\ot q_{j_{1}l_{1}}...q_{j_{\gamma+1}l_{\gamma+1}})(\sum_{(m_{1},...,m_{2\gamma+2}) \ \text{path}}(q_{i_{1}m_{1}}...q_{i_{\gamma+1}m_{\gamma+1}}...q_{i_{2\gamma+2}m_{2\gamma+2}}\ot\\
		&& q_{m_{1}l_{1}}...q_{m_{\gamma+1}l_{\gamma+1}}...q_{m_{2\gamma+2}l_{2\gamma+2}})(q_{i_{\gamma+2}j_{\gamma+2}}...q_{i_{2\gamma+2}j_{2\gamma+2}}\ot q_{j_{\gamma+2}l_{\gamma+2}}...q_{j_{2\gamma+2}l_{2\gamma+2}})\\
		&=&\sum_{(m_{1},...,m_{2\gamma+2}) \ \text{path}}( (q_{i_{1}j_{1}}...q_{i_{\gamma+1}j_{\gamma+1}}q_{i_{1}m_{1}}...q_{i_{\gamma+1}m_{\gamma+1}})(q_{i_{\gamma+2}m_{\gamma+2}}...q_{i_{2\gamma+2}m_{2\gamma+2}}q_{i_{\gamma+2}j_{\gamma+2}}...\\
		&&q_{i_{2\gamma+2}j_{2\gamma+2}})\ot
		(q_{j_{1}l_{1}}...q_{j_{\gamma+1}l_{\gamma+1}}q_{m_{1}l_{1}}...q_{m_{\gamma+1}l_{\gamma+1}})(q_{m_{\gamma+2}l_{\gamma+2}}...q_{m_{2\gamma+2}l_{2\gamma+2}} q_{j_{\gamma+2}l_{\gamma+2}}\\
		&&...q_{j_{2\gamma+2}l_{2\gamma+2}}))
		\end{eqnarray*} 
	By Lemma \ref{key} and (\ref{1}), the last summation reduces to
	\begin{displaymath}
	q_{i_{1}j_{1}}...q_{i_{\gamma+1}j_{\gamma+1}}q_{i_{\gamma+2}j_{\gamma+2}}...q_{i_{2\gamma+2}j_{2\gamma+2}}\ot q_{j_{1}l_{1}}...q_{j_{\gamma+1}l_{\gamma+1}} q_{j_{\gamma+2}l_{\gamma+2}}...q_{j_{2\gamma+2}l_{2\gamma+2}}.
	\end{displaymath}
	Hence
	\begin{eqnarray*}
	&&(q_{i_{1}j_{1}}...q_{i_{\gamma+1}j_{\gamma+1}}\ot q_{j_{1}l_{1}}...q_{j_{\gamma+1}l_{\gamma+1}})\Delta(q_{i_{1}l_{1}}...q_{i_{\gamma+1}l_{\gamma+1}}...q_{i_{2\gamma+2}l_{2\gamma+2}})\\
	&&(q_{i_{\gamma+2}j_{\gamma+2}}...q_{i_{2\gamma+2}j_{2\gamma+2}}\ot q_{j_{\gamma+2}l_{\gamma+2}}...q_{j_{2\gamma+2}l_{2\gamma+2}})\\
	&=& q_{i_{1}j_{1}}...q_{i_{\gamma+1}j_{\gamma+1}}q_{i_{\gamma+2}j_{\gamma+2}}...q_{i_{2\gamma+2}j_{2\gamma+2}}\ot q_{j_{1}l_{1}}...q_{j_{\gamma+1}l_{\gamma+1}} q_{j_{\gamma+2}l_{\gamma+2}}...q_{j_{2\gamma+2}l_{2\gamma+2}}\\
	&\neq& 0, 
	\end{eqnarray*}
since by assumption, $(i_{1},...,i_{2\gamma+2})\sim^{k}(j_{1},...,j_{2\gamma+2})\sim^{k}(l_{1},...,l_{2\gamma+2})$. So \begin{displaymath}q_{i_{1}l_{1}}...q_{i_{\gamma+1}l_{\gamma+1}}...q_{i_{2\gamma+2}l_{2\gamma+2}}\neq 0,\end{displaymath} proving the required transitivity. 
\end{proof}
We already know that for a fixed $k$, $\mathbb{G}^{+}_{k}(\Gamma)$ acts on $C(E^{\gamma})$ for all $\gamma\leq k$. Also we know that there is a $\sim^{k}$ equivalence relation on $E^{\gamma}$. This allows us to state and prove the following lemma.
\blmma
\label{erg1}
For any $\gamma\leq k$, consider the coaction $u^{\ot(\gamma+1)}:C(E^{\gamma})\raro C(E^{\gamma})\ot C(\mathbb{G}^{+}_{k}(\Gamma))$. Then $u^{\ot(\gamma+1)}(f)=f\ot 1$ if and only if $f$ is constant on the $\sim^{k}$ orbits.
\elmma
\begin{proof}
	The proof is given along the same lines of arguments used in \cite{Maccinska}. Let us write $f$ as $\sum_{(i_{1},...,i_{\gamma+1}) \ {\rm path}}f^{i_{1}...i_{\gamma+1}}\chi_{i_{1}...i_{\gamma+1}}$. Then
	\begin{displaymath}
		u^{\ot(\gamma+1)}(f)=\sum_{(j_{1},...,j_{\gamma+1}) \ {\rm path}}\chi_{j_{1}...j_{\gamma+1}}\ot\sum_{(i_{1},...,i_{\gamma+1}) \ path}f^{i_{1},...,i_{\gamma+1}}q_{j_{1}i_{1}}...q_{j_{\gamma+1}i_{\gamma+1}}.
		\end{displaymath}
	Since $u^{\ot(\gamma+1)}(f)=f\ot 1$, comparing coefficints, we get
	\begin{displaymath}
	f^{j_{1}...j_{\gamma+1}}=\sum_{(i_{1}...i_{\gamma+1}) {\rm path}}f^{i_{1}...i_{\gamma+1}}q_{j_{1}i_{1}}...q_{j_{\gamma+1}i_{\gamma+1}} \forall \ (j_{1},...,j_{\gamma+1}) \ {\rm path}.
	\end{displaymath}
	For fixed paths $(j_{1},...,j_{\gamma+1})$ and $(l_{1},...,l_{\gamma+1})$, multiplying both sides with $q_{j_{1}l_{1}}...q_{j_{\gamma+1}l_{\gamma+1}}$ and using Lemma \ref{key}, we get
	\begin{displaymath}
	(f^{j_{1}...j_{\gamma+1}}-f^{l_{1}...l_{\gamma+1}})q_{j_{1}l_{1}}...q_{j_{\gamma+1}l_{\gamma+1}}=0.
	\end{displaymath}
	If $(j_{1},...,j_{\gamma+1})\sim^{k}(l_{1},...,l_{\gamma+1})$, then $f^{j_{1}...j_{\gamma+1}}=f^{l_{1}...l_{\gamma+1}}$, proving the claim.
\end{proof}
For notations and proof of the following corollary the reader is referred to \cite{bancolo}.
\bcrlre
For $k\in\mathbb{N}$, fix any $\gamma\leq k$. Then TFAE\\
(i) The action $u^{\ot(\gamma+1)}$ of $\mathbb{G}^{+}_{k}$ on $C(E^{\gamma})$ is $\sim^{k}$ transitive.\\
(ii) The action is ergodic i.e. $\mathrm{Fix}(u^{\ot(\gamma+1)})=\mathbb{C}1$.\\
(iii) $\sum_{(i_{1},...,i_{\gamma+1}) \ {\rm path}}\int_{\mathbb{G}^{+}_{k}(\Gamma)}q_{i_{1}i_{1}}...q_{i_{\gamma+1}i_{\gamma+1}}=1$.
\ecrlre
\subsection{An application: Bound on the non zero eigenvalues of the graph Laplacian}
In this subsection we shall improve the bound of the first non zero eigenvalue of $\cll$ of a finite, connected vertex transitive graph without multiple edges and loops. We shall adapt the techniques used in \cite{Chung} to the quantum set up to improve the lower bound (\ref{bound}). To that end, given a graph $\Gamma$, recall the equivalence relation $\sim^{k}$ on the edge set. Let $\{E_{i}^{k}:i=1,...,m\}$ be the set of $\sim^{k}$-equivalence classes of edges. We define $r_{k}$ to be $\mathrm{min}_{i}|E_{i}^{k}|$. Then $r_{k}\geq r_{k+1}$ for all $k$. If we define $\mathrm{ind}_{k}(\Gamma)$ by $\frac{\clv}{r_{k}}$, we immediately have
\begin{eqnarray} \label{chain2}\mathrm{ind}_{k}(\Gamma)\leq \mathrm{ind}_{k+1}(\Gamma)\leq...\leq{\rm ind}(\Gamma) \ for \ all \ k.\end{eqnarray}\\
\indent Let $\Gamma$ be a finite, connected, vertex transitive graph of diameter $D\leq (2k+1)$ and $n$ number of vertices. Fix a vertex $i_{1}\in V$ and consider the following set
\begin{displaymath}
\cls=\{(i_{1},...,i_{\alpha+1}): (i_{1},...,i_{\alpha+1}) \ \mathrm{is \ a \ fixed \ shortest \ path}\}.
\end{displaymath}
Note that $\alpha$ is not a fixed integer, but the assumption on the diameter of the graph ensures that $\alpha\leq(2k+1)$. Now for a vertex $j_{1}\in V$ let us define the following set
\begin{displaymath}
\clp_{k}^{j_{1}}=\{(j_{1},...,j_{\alpha+1}):(j_{1},...,j_{\alpha+1})\sim^{k}(i_{1},...,i_{\alpha+1}), (i_{1},...,i_{\alpha+1})\in\cls\}.
\end{displaymath}
The vertex transitivity ensures that $\clp_{k}^{j_{1}}$ is non empty for all $j_{1}\in V$. Finally for an edge $e$, generalizing the definition of $N_{e}$ in \cite{Chung}, we define
\begin{displaymath}
N^{k}_{e}=\#\{(j_{1},...,j_{\alpha+1})\in\clp_{k}^{j_{1}}:j_{1}\in V, e \ \mathrm{occurs \ in} \ (j_{1},...,j_{\alpha+1})\}.
\end{displaymath}
Using the fact that ${\rm Aut}(\Gamma)$ is a quantum subgroup of $\mathbb{G}^{+}_{k}$ for all $k\in\mathbb{N}$, it is easy to see that $N^{k}_{e}\geq N_{e}$ for all $k$.
\bppsn
\label{same length}
For two edges $e\sim^{k}f$, we have $N_{e}^{k}=N_{f}^{k}$.
\eppsn
\begin{proof} 
Let $e=(j_{i},j_{i+1})\sim^{k}(l_{i},l_{i+1})=f$. We shall show that $e$ occurs in some path in $\clp_{k}^{j_{1}}$ for some $j_{1}\in V$ if and only if $f$ occurs in some path in $\clp_{k}^{l_{1}}$ for some $l_{1}\in V$ proving the proposition. To that end let $(j_{1},...j_{i},j_{i+1},...,j_{\alpha+1})\in \clp_{k}^{j_{1}}$. We claim that there exists a path $(l_{1},...,l_{i},l_{i+1},...,l_{\alpha+1})$ such that $(j_{1},...j_{i},j_{i+1},...,j_{\alpha+1})\sim^{k}(l_{1},...l_{i},l_{i+1},...,l_{\alpha+1})$. Then since $\alpha\leq(2k+1)$, by transitivity of $\sim^{k}$, we can conclude in deed $f$ occurs in $(l_{1},...l_{i},l_{i+1},...,l_{\alpha+1})\in\clp_{k}^{l_{1}}$. Let us prove the claim. Since $(j_{i},j_{i+1})\sim^{k}(l_{i},l_{i+1})$, we have $q_{j_{i}l_{i}}q_{j_{i+1}l_{i+1}}\neq 0$ for $q_{j_{i}l_{i}},q_{j_{i+1}l_{i+1}}\in C(\mathbb{G}^{+}_{k}(\Gamma))$. If possible, suppose $q_{j_{1}l_{1}^{\prime}}...q_{j_{i}l_{i}}q_{j_{i+1}l_{i+1}}...q_{j_{\alpha+1}l_{\alpha+1}^{\prime}}=0$ for all paths $(l_{1}^{\prime},...,l_{i},l_{i+1},...,l_{\alpha+1}^{\prime})$. Hence
\begin{displaymath}
\sum_{(l_{1}^{\prime},...,l_{\alpha+1}^{\prime}) \ {\rm path}}q_{j_{1}l_{1}^{\prime}}...q_{j_{i}l_{i}}q_{j_{i+1}l_{i+1}}...q_{j_{\alpha+1}l_{\alpha+1}^{\prime}}=0.
\end{displaymath}
Again we can replace the path indices by arbitrary indices as before and hence we get using the magic unitarity, $q_{j_{i}l_{i}}q_{j_{i+1}l_{i+1}}=0$, a contradiction. So in deed there is a path $(l_{1}^{\prime},...,l_{i},l_{i+1},...,l_{\alpha+1}^{\prime})$ with $q_{j_{1}l_{1}^{\prime}}...q_{j_{i}l_{i}}q_{j_{i+1}l_{i+1}}...q_{j_{\alpha+1}l_{\alpha+1}^{\prime}}\neq 0$, completing the proof of the proposition.
\end{proof}
So we can estimate $N^{k}_{e}$ along the lines of the proof of the Theorem 7.7 of \cite{Chung}.
\begin{eqnarray}
\label{auxbound}
N^{k}_{e}\leq\frac{n^{2}D}{2|E_{i}^{k}|}\leq\frac{n^{2}D}{2\mathrm{min}|E_{i}^{k}|}\leq\frac{nD\mathrm{ind}_{k}(\Gamma)}{\clv}.
\end{eqnarray}
Now we are ready to prove the main result of this subsection.
\bthm
\label{main}
Let $\Gamma$ be a finite, connected vertex transitive graph without multiple edge or loops of diameter $D\leq(2k+1)$. Then we have the following lower bound of the first non-zero eigenvalue of $\Gamma$:
\begin{eqnarray}
\label{mainres}
\lambda_{1}\geq \frac{1}{D^{2}\mathrm{ind}_{k}(\Gamma)}.
\end{eqnarray}
\ethm 
\begin{proof}
	Once again, we shall adapt the techniques used in Theorem 7.7 of \cite{Chung}. Since the graph is vertex transitive, we assume the graph to be $s$-regular. Then $s=\frac{\clv}{n}$. Let us write the first non zero eigenvalue as (see Theorem 7.7 of \cite{Chung})
	\begin{eqnarray*}
		\lambda_{1}={\rm min}_{f}\frac{n\sum_{i\sim j}(f(i)-f(j))^{2}}{s\sum_{i,j}(f(i)-f(j))^{2}}
		\end{eqnarray*}
	Again as in Theorem 7.7 of \cite{Chung}, writing $f(e)=|f(i)-f(j)|$ for $e=(i,j)$, we can rewrite $\lambda_{1}$ as
	\begin{displaymath}
	{\rm min}_{f}\frac{n\sum_{e\in E}f^{2}(e)}{s\sum_{i,j}(f(i)-f(j))^{2}}.
	\end{displaymath}
	Now fix a harmonic eigenfunction $f$ achieving the minimum. We have (see page 118 of \cite{Chung}), $\sum_{i,j}(f(i)-f(j))^{2}\leq f^{2}(e)DN_{e}$. Since $N^{k}_{e}\geq N_{e}$, we have
	\begin{displaymath}
	\sum_{i,j}(f(i)-f(j))^{2}\leq\sum_{e\in E} f^{2}(e)DN_{e}^{k}\leq \sum_{e\in E} f^{2}(e)D\frac{D{\rm ind}_{k}(\Gamma)}{s} (by \ (\ref{auxbound})).
	\end{displaymath}
	Therefore we have
	\begin{eqnarray*}
		\lambda_{1}&=&\frac{n\sum_{e\in E}f^{2}(e)}{s\sum_{i,j}(f(i)-f(j))^{2}}\\
		&\geq& \frac{1}{D^{2}\mathrm{ind}_{k}(\Gamma)}.
		\end{eqnarray*}
\end{proof}
\brmrk
1. From (\ref{chain2}) and (\ref{mainres}) we see that the lower bound of the first non zero eigenvalue gets stronger as $k$ decreases. But with decreasing $k$, the applicabilty of Theorem \ref{main} gets restricted because of the restriction on the diameter of the graph.\\
\indent 2. It would be interesting to find concrete examples of graphs where we can see that Theorem \ref{main} indeed gives better bound on the lower bound of the first non zero eigenvalue. The techniques used in \cite{Maccinska} to find example of quantum vertex transitive, but not vertex transitive graph might be useful in this regard. The author hopes that this quest would also open up the study of quantum edge transitive graphs (in a suitable sense) which are not edge transitive in the classical sense.
\ermrk
{\bf Acknowledgement}: The author would like to thank both Arnab Mandal and Anirban Banerjee for some fruitful discussions.


\begin{thebibliography} {andrw}
	\bibitem{banica} T. Banica, {\it Quantum automorphism groups of homogeneous graphs},  J. Funct. Anal.,  224(2005), 243-280.
	\bibitem{bancolo} T. Banica, {\it Higher transitive quantum groups: theory and models}, Colloq. Math., 156(2019), no. 1, 1-14.
	\bibitem{bichon} J. Bichon: {\it Quantum automorphism groups of finite graphs},  Proc. Amer. Math. Soc.,  131(2003), no 3, 665-673.
	\bibitem{Chung} R. K. Chung: {\it Spectral graph theory}, CBMS, number 92.
	\bibitem{Maccinska} M. Lupini, L. Mancinska, David E. Roberon,: {\it Nonlocal games and quantum permutation groups}, arXiv: 1712.01820.
	\bibitem{Van} A. Maes A. and  A. Van daele, {\it Notes on compact quantum groups},  Nieuw Arch. Wisk (4)  16 (1998), no.1-2, 73-112 .
	\bibitem{Schmidt} S. Schmidt, M. Weber, {\it Quantum symmetries of graph $C^{\ast}$-algebras}, Canad. Math. Bull., 61(2018), no 4, 848-864.
	\bibitem{wang} S. Wang, {\it Quantum symmetry groups of finite spaces},  Comm. Math. Phys.,  195(1998), 195-211.
	\bibitem{Woro} S.L. Woronowicz, {\it Compact matrix pseudogroups},  Comm. Math. Phys.,  111(1987), 613-665.
	
\end{thebibliography}
\end{document}